%cd Desktop
%amstex h2.tex
%dvipdfm h2.dvi
%tex h2
%dvips h2
%ps2pdf h2.ps
\input amstex
\documentstyle{amsppt}
\magnification1200
\tolerance=10000
\def\n#1{\Bbb #1}
\def\p{\Bbb C_{\infty}}

\def\Tr{\hbox{Tr}}

\def\End{\hbox{End }}

\def\Ker{\hbox{Ker }}

\def\Res{\hbox{Res }}

\def\rd{\sqrt{-\Delta}}

\def\e11{E_{11}}

\def\vf{\varphi}

\topmatter
\title
Anderson T-motives and abelian varieties with MIQF: results coming from an analogy
\endtitle
\author
A. Grishkov, D. Logachev\footnotemark \footnotetext{E-mails: shuragri{\@}gmail.com; logachev94{\@}gmail.com (corresponding author)\phantom{**************}}
\endauthor
\NoRunningHeads
\address
First author: Departamento de Matem\'atica e estatistica
Universidade de S\~ao Paulo. Rua de Mat\~ao 1010, CEP 05508-090, S\~ao Paulo, Brasil, and Omsk State University n.a. F.M.Dostoevskii. Pr. Mira 55-A, Omsk 644077, Russia.
\medskip
Second author: Departamento de Matem\'atica, Universidade Federal do Amazonas, Manaus, Brasil
\endaddress
\thanks Thanks: Alain Genestier and Laurent Fargues attracted the attention of the authors to
the main analogy of the present paper. The authors are grateful to FAPESP, S\~ao Paulo, Brazil for a financial support (process No. 2017/19777-6). The first author is grateful to SNPq, Brazil, to RFBR, Russia, grant 16-01-00577a (Secs. 1-4), and to Russian Science Foundation, project 16-11-10002 (Secs. 5-8) for a financial support. The authors are grateful to an anonymous reviewer for some important remarks.\endthanks

\abstract Analogy between Anderson T-motives and abelian varieties with multiplication by an imaginary quadratic field (MIQF) is a source of 2 results:

1. A description of abelian varieties with MIQF of dimension $r$ and signature $(n, r-n)$ in terms of  "lattices" of dimension $r$ in $\n C^n$;

2. A construction of exterior powers of abelian varieties with MIQF having $n=1$.
\endabstract
\endtopmatter
\document
{\bf 0. Introduction.}
\nopagebreak
\medskip
The origin of the present paper is an analogy 1.8 (known for experts) between the following two objects
corresponding to the function field case and the number field case respectively:
\medskip
{\bf A.} A T-motive $M$ of rank $r$ and dimension $n$, pure,
uniformizable, having the nilpotent operator $N$ equal to 0 (see Section 1 for the exact
definitions and the origin of the analogy),
\medskip
and
\medskip
{\bf B.} An abelian variety $A$ over $\n C$ with multiplication by an imaginary quadratic
field $K$ (abbreviated as an abelian variety with MIQF), of dimension $r$ and of
signature $(n, r-n)$.
\medskip
Using well-known constructions for T-motives, we get analogous
constructions for abelian varieties with MIQF. Firstly, let us consider the lattice
$L(M)$ associated to a T-motive $M$ of type {\bf A}:
\medskip
{\bf C.} $L(M)$ is an $r$-dimensional $\n F_p[\theta]$-lattice in $\p^n$, where $\n F_p [\theta]$ (resp. $\p$) is the functional analog of $\n Z$ (resp. $\n C$).
\medskip
Let us consider the lattice $L(A)$ associated to an abelian variety $A$ with MIQF:
\medskip
{\bf D.} $L(A)$ is an $r$-dimensional $O_K$-lattice in $\n C^r$ ($L(A)$ is not an
$O_K$-submodule of $\n C^r$ treated as $O_K$-module).
\medskip
Objects of types {\bf (C)} and {\bf (D)} are not quite analogous, isn't it? The first
result of the present paper is Theorem 2.6 which gives us --- roughly speaking --- that
\medskip
{\bf E.} $A$ defines an $r$-dimensional $O_K$-submodule
of $\n C^n$ (not of $\n C^r$ ! ).
\medskip
{\bf (E)} is an analog of {\bf (C)}.
\medskip
Secondly, it is known that if $M$ of type {\bf A} has $n=1$ then its $k$-th exterior
power $\lambda^k(M)$ is also an object of type {\bf A}. By analogy, we can expect that if
$A$ is an abelian variety with MIQF having $n=1$ then $\lambda^k(A)$ is defined, and is
also an abelian variety with MIQF. We really give this construction (Section 3).
\medskip
Both these results are of elementary nature, they could be known to Riemann.\footnotemark \footnotetext{ Claire Voisin told me that she knew the construction of Section 3 (not published), but she did not know its relation with exterior power of Drinfeld modules.}
Probably this analogy will be a source of more results. For example, it would be interesting to
define analogs of T-motives having $N\ne0$.
\medskip
The paper is organized as follows. In Section 1 we recall definitions of
T-motives and we formulate the results which are starting points (using the analogy 1.8)
of the results of Sections 2, 3.  Section 2 contains the exact statement and the proof
of {\bf (E)}. In Section 3 we apply this result to
construct exterior powers of abelian varieties with MIQF having $n=1$. Formally, these
sections are independent of Section 1, i.e. they do not require any knowledge of
functional case.
\medskip
{\bf Remark.} Many constructions/results of the theory of Drinfeld modules are analogs of the corresponding constructions/results of the theory of abelian varieties. The direction of analogies of the present paper is the opposite: from Drinfeld modules to abelian varieties.
\medskip
{\bf Section 1. Origin of construction: T-motives.}
\medskip
A standard reference for T-motives is [G], Section 5, we shall use its notations if
possible.
Let $\goth r$ be a power of a prime number, \footnotemark \footnotetext{We use notation
$\goth r$ instead of
$r$ of [G] in order do not confuse with the rank of a T-motive.} $\n F_{\goth r}(\theta)$
(resp. $\n F_{\goth r}((\theta^{-1}))$) the functional analog of $\n Q$ (resp. of $\n
R$), and
$\p$ --- the completion of the algebraic closure of $\n F_{\goth r}((\theta^{-1}))$ ---
the functional
analog of $\n C$.

The definition of a T-motive $M$ is given in [G], Definitions 5.4.2,
5.4.12 (Goss uses another terminology: "abelian T-motive" of [G] = "T-motive" of the present paper; $L$ of Goss should be considered as $\p$). Particularly, $M$ is a free $\p[T]$-module of
dimension $r$ (this number $r$ is called the rank of $M$) endowed by a $\p$-skew-linear
operator $\tau$ satisfying some properties. A nilpotent operator $N=N(M)$
associated to a T-motive is defined in [G], Remark 5.4.3.2. We shall consider only pure ([G],
Definition 5.5.2) uniformizable ([G], Theorem 5.9.14, (3)) T-motives. Its dimension $n$
is defined in [G], Remark 5.4.13.2 (Goss denotes the dimension by $\rho$). Condition
$N=0$ implies $n \le r$. A T-motive of dimension 1 is the same as a
Drinfeld module, they are all pure, uniformizable, and their $N$ is 0.

If $N(M)=0$ then attached to such T-motive
is a lattice $L=L(M)$ which is a free $r$-dimensional $\n F_{\goth r}[\theta]$-module in
$\p^n$, and if $N\ne 0$ then $L(M)$ is a slightly more complicated object, we do not need
to consider details for this case.
Inclusion of $L$ in $\p^n$ defines a surjective map
$$\alpha=\alpha(M): L\underset{\n F_{\goth r}[\theta]} \to{\otimes}\p \to \p^n
\eqno{(1.1)}$$
having the property
$$\hbox{ Restriction of $\alpha$ to $L\underset{\n F_{\goth r}[\theta]} \to{\otimes} \n
F_{\goth r}((\theta^{-1}))$
is injective} \eqno{(1.2)}$$

Tensor product of T-motives $M_1$, $M_2$ is simply their tensor product
over $\p[T]$, where the action of $\tau$ is defined by the formula $\tau(m_1 \otimes
m_2)=\tau(m_1) \otimes \tau(m_2)$. If both $M_1$, $M_2$ are pure uniformizable then
$M_1\otimes M_2$ is pure uniformizable as well. The same definition holds for exterior
(resp. symmetric) powers of $M$. The dual T-motive $M^*$ is defined in [GL07].
\medskip
{\bf 1.3.} It is easy to check that even if $N(M_1)$, $N(M_2)$ are 0 then $N(M_1\otimes
M_2)$, $N(\lambda^k(M_1))$, $N(S^k(M_1))$ are not 0. The only exception: if $M$ is a
Drinfeld module ($\iff n=1$) then $N(\lambda^k(M))=0$, this is an elementary calculation.
\medskip
There is a natural problem to describe $L(M_1\otimes M_2)$, $L(\lambda^k(M_1))$,
$L(S^k(M_1))$ in terms of $L(M_1)$, $L(M_2)$. It was solved by Anderson (non-published),
a formula which is equivalent to this description is stated without proof in [P], end of
page 3. See [GL07], Section 6 for more details and for the proof of the formula of Pink
in the case when $N(M_1)$, $N(M_2)$ are 0. Let us state the theorem for the case when all
$N$ are 0 (and hence all $L$ are lattices), i.e. $M$ is a Drinfeld module of rank $r$. We
need
\medskip
{\bf Definition 1.4.} For a short exact sequence of vector spaces over a field
$$0\to B_1\overset{i}\to{\to } B_2\to C\to
0$$ we define its
$k$-th exterior power as the following exact sequence:
$$0\to \lambda^k(B_1)\overset{\lambda^k(i)}\to{\to } \lambda^k(B_2)\to C_k\to 0$$

Now let us consider the exact sequences for $M$, $\lambda^k(M)$
$$0\to \Ker \alpha(M) \to L(M)\underset{\n F_{\goth r}[\theta]} \to{\otimes}\p
\overset{\alpha(M)}\to{\to}  \p\to 0\eqno{(1.5)}$$
$$0\to \Ker \alpha(\lambda^k(M)) \to L(\lambda^k(M))\underset{\n F_{\goth r}[\theta]}
\to{\otimes}\p \overset{\alpha(\lambda^k(M))}\to{\to} \p^{n(r,k)} \to 0\eqno{(1.6)}$$
where maps $\alpha(M)$, $\alpha(\lambda^k(M))$ are from 1.1, and $n(r,k)= {r-1 \choose k-1}$ is the dimension of $\lambda^k(M)$.
\medskip
{\bf Theorem 1.7.} There exists a canonical isomorphism from the $k$-th exterior power of
(1.5) to (1.6) such that the image of $\lambda^k(L(M))$ in (1.5) is $L(\lambda^k(M))$ in
(1.6).
\medskip
{\bf Proof} is completely analogous to the proofs for the case of dual
T-motives ([GL07], Theorem 5) and for the case of tensor product of
T-motives having $N=0$ ([GL07], Theorem 6.2), hence it is omitted. $\square$
\medskip
{\bf 1.8. Origin of the analogy.} We give here only a sketch of definitions; see for
example [W] for the exact statements. Let $X$ be a Shimura variety and $G$ a reductive
group over $\n Q$ associated to $X$ according Deligne. Let $p$ be a prime of good
reduction of $X$. $\n H_p(X)$ --- the $p$-part of the Hecke algebra of $X$ --- is
isomorphic to $\Cal H(G(\n Q_p))$ --- the algebra of double cosets $G(\n Z_p)gG(\n Z_p)$,
$g\in G(\n Q_p)$. There exists a Levi subgroup $M$ of $G$ having the following property:
\medskip
The $p$-part of the Hecke algebra of $\tilde X$ (the reduction of $X$ at $p$) is
isomorphic to $\Cal H(M(\n Q_p))$.
\medskip
See, for example, [W], p. 44, (*) and p. 49, (1.10) for the definition and properties of
$M$.
\medskip
{\bf Example.} If $X$ is a Shimura variety parametrizing abelian varieties with
multiplication by an imaginary quadratic field $K$, of dimension $r$ and of signature
$(n, r-n)$, then $G(\n Q)=GU(n,r-n)(K)$. If $p$ splits in $K/\n Q$ then $G(\n
Q_p)=GL_r(\n Q_p)\times GL_1(\n Q_p)$. Let $i_{n,r-n}: GL_n(\n Q_p)\times GL_{r-n}(\n Q_p) \hookrightarrow GL_r(\n Q_p)$ be the inclusion
of the subgroup of $(n\times r-n)$-block diagonal matrices. Then the inclusion $M(\n Q_p) \hookrightarrow G(\n Q_p)$ is the inclusion
$$i_{n,r-n}\times GL_1(\n Q_p): GL_n(\n Q_p)\times GL_{r-n}(\n Q_p)\times GL_1(\n Q_p) \hookrightarrow GL_r(\n Q_p)\times GL_1(\n Q_p)$$

For the function field case an analog of this theory is conjectural, but preliminary results
of [L] show that for abelian Anderson T-motives of rank $r$ and dimension $n$ we have
the same groups (up to the factor $GL_1$): $G=GL_r$, $M= GL_n \times GL_{r-n}\subset GL_r$.
\medskip
As a corollary we get that the dimensions of moduli spaces of both types of objects
(T-motives; abelian varieties with MIQF) are equal: they are $n(r-n)$.
\medskip
The below sections 2, 3 contain constructions of the number field case analogs of the map
$\alpha$ of (1.1), and of Theorem 1.7 respectively. From one side, finding of these
constructions was inspired by the analogy; from another side, their existence is a
support to the analogy.
\medskip
{\bf 2. Abelian varieties with MIQF.} We shall fix an imaginary quadratic field $K=\n
Q(\rd)$. For simplicity, an abelian variety $A$ is treated up to isogeny, and we restrict
ourselves only by one fixed polarization form.

The main theorem 2.6 establishes an equivalence of objects of types {\bf B} and {\bf E}.
\medskip
{\bf Definitions for the type B.} Let $A=V/D_{\n Z}$, $V=\n C^r$ be an abelian
variety with MIQF,
$L=D_{\n Z}\otimes_{\n Z}{\n Q}$. Since we consider $A$ up to isogeny, we shall deal only
with $L$ and not with
$D_{\n Z}$. We fix an inclusion $\iota: K\hookrightarrow \End(A)\otimes_{\n Z} \n Q$
defining multiplication, and we fix an Hermitian polarization form $H=B+i\Omega$ of $A$
on $V$,
where $B$ and $\Omega$ are respectively its real and imaginary parts. There are two
structures of $K$-module
on $V$: the ordinary one which is the restriction of the $\n C$-module structure, and the
*-structure
(multiplication is denoted by $k*v$, $k\in K$, $v \in V$) coming from $\iota$. $L$ is a
$K$-*-module. We choose a basis $\goth x_1,...,\goth x_r$ of $L/K$ (notations of [Sh]).
According [Sh], p. 157,
(11) there exists a matrix $T=\{t_{ij}\}\in M_r(K)$ such that
$$\forall k_1, k_2\in K, i,j=1,...,r \ \ \Omega(k_1*\goth x_i, k_2*\goth x_j)=\Tr_{K/\n
Q}(k_1t_{ij}\bar k_2)
\eqno{(2.1)}$$
([Sh], p. 157, (11)). $T$ has properties
\medskip
(a) $\bar T^t=-T$, i.e. $iT$ is hermitian ([Sh], p. 157, (12)) and
\medskip
(b) Signature of $iT$ is $(n, r-n)$ ([Sh], p. 160, (25)).
\medskip
We restrict ourselves by those $A$ whose $T$ (it depends on $\goth x_1,...,\goth x_r$)
satisfy
$$T=\rd E_{n,r-n}\eqno{(2.2)}$$ where $E_{n,r-n}:=\left(\matrix E_n & 0 \\
0&-E_{r-n}\endmatrix \right)$.
\medskip
{\bf 2.3. Definitions for the type E.} We consider the set of triples $(L, H_L,
\alpha)$ where
\medskip
(1) $L$ is a $K$-vector space of dimension $r$;
\medskip
(2) $H_L$ is a $K$-valued Hermitian form on $L$ of signature $(n, r-n)$ such that there
exists a basis of $L$ over $K$ satisfying the condition:
$$\hbox{the matrix of $H_L$ in this basis is $E_{n,r-n}$}\eqno{(2.4)}$$

\medskip
We denote by $H_{L, \n C}$ the extension of $H_L$ to $L\otimes_K\n C$.
\medskip
(3) $\alpha: L\otimes_K\n C\to \n C^n$ is a $\n C$-linear map such that $\alpha$ is
surjective, and
$$\hbox{The restriction of $-H_{L,\n C}$ to $\Ker \alpha$ is a positive definite
form.}\eqno{(2.5)}$$

{\bf Remark.} This $\alpha$ is clearly an analog of $\alpha$ of 1.1. We see that
conditions of surjectivity of $\alpha$ hold in both cases, while the property 1.2
apparently has no analog in the number field case.
\medskip
{\bf Theorem 2.6.} There is a 1 -- 1 correspondence between the above $A$, $\iota$, $H$
--- objects of type {\bf (B)} (here $A$ is up to
isogeny, and $A$ has $T$ satisfying
2.2), and the above triples $(L, H_L, \alpha)$ --- objects of type {\bf (E)}.
\medskip
{\bf Remark.} If the restriction of $\alpha$ on $L \subset L\otimes_K\n C$ is injective
(this holds for almost all triples $(L, H_L, \alpha)$ ) then the $K$-vector space in $\n
C^n$ mentioned in {\bf (E)} is $\alpha(L)$.
\medskip
{\bf Proof.} $L$ is the same for both types (when we consider $L$ for the type {\bf (E)},
we omit * in $k*l$). Let the triple ($A$, $\iota$, $H$) of the type {\bf (B)} be given.
There is a canonical decomposition $V=V^+\oplus V^-$ where
$$V^+:=\{v\in V \ \vert \ k*v=kv\}\eqno{(2.7)}$$
$$V^-:=\{v\in V \ \vert \ k*v=\bar kv\}\eqno{(2.8)}$$

We denote by $\pi_+:V\to V^+$ the projection along $V^-$. Let $\beta: L \to V$ be the
tautological inclusion.
\medskip
{\bf 2.9.} The triple $(L, H_L, \alpha)$ corresponding to ($A$, $\iota$, $H$) is
constructed as
follows. We define
$$\alpha:L\otimes_K\n C\to V^+=\n C^n$$ by the formula $\alpha(l\otimes
z)=z\cdot\pi_+(\beta(l))$. Formula 2.7 shows that it is well-defined. The form $H_L$ is
defined by the equality
$$\forall l_1, l_2\in L \ \ \ \Tr_{K/\n Q}(\rd H_L(l_1,l_2))=\Omega(\beta(l_1),
\beta(l_2))\eqno{(2.10)}$$

{\bf Lemma.} Formula 2.10 really defines $H_L$ uniquely.
\medskip
{\bf Proof. Unicity:} we fix $l_1, l_2\in L$, and we consider a $\n Q$-linear form
$\gamma: K \to \n Q$ defined by
the formula

$$\gamma(k) = \Omega(\beta(k*l_1), \beta(l_2))$$
There exists the only $k_0\in K$ such that $\gamma(k) = \Tr_{K/\n Q}(k_0k)$. We see that
if 2.10 holds for all pairs
$k*l_1, l_2$, then necessarily $\rd H_L(l_1,l_2)=k_0$.
\medskip
{\bf Existence:} Let $\goth x_1,...,\goth x_r$ be a basis of $L/K$ such that the
corresponding
$T$ has the form 2.2.
We define $H_L$ by the condition that the matrix of $H_L$ in $\goth x_1,...,\goth x_r$ is
$E_{n,r-n}$. 2.1
implies that $H_L$ satisfies 2.10. $\square$
\medskip
To prove 2.5 we recall a well-known
\medskip
{\bf 2.11. Coordinate description.} The Siegel domain
for the present case is the following:
$$\hbox{$\Cal H^3_{n,r-n}:=\{z\in M_{n,r-n}(\n C)\ \vert \ E_n-z\bar z^t$ is positive
hermitian $\}$ }\eqno{(2.12)}$$
$$( \hbox{$\iff \ E_{r-n}-\bar z^t z$ is positive
hermitian })$$
([Sh], p. 162, 2.6). For any $z\in \Cal H^3_{n,r-n}$ we can construct an above abelian
variety $A_z$ (satisfying 2.2) as follows ([Sh]). We fix a basis $e_1,...,e_r$ of $V$
over $\n C$
such that $e_1,...,e_n$ (resp. $e_{n+1},...,e_r$) is a basis of $V^+$ (resp. $V^-$) over
$\n C$. It
defines the $K$-*-action on $V$. Let $\goth x_*$ (resp. $e_*$) be the matrix column
of $\goth x_1,...,\goth x_r$ (resp.
$e_1,...,e_r$). They satisfy
$$\goth x_*=Y e_* \hbox{ where } Y=\left(\matrix E_n & z \\
\bar z^t&E_{r-n}\endmatrix \right) \eqno{(2.13)}$$ ([Sh], p. 162, (35), Type IV). $\Omega$
is defined by 2.1,
2.2. These conditions define $A_z$.
\medskip
{\bf Proof of 2.5.} 2.13 implies that elements
$$\lambda_i:=\goth x_{n+i}\otimes 1 -\sum_{k=1}^n \goth x_k\otimes\bar  z_{ki}\eqno{(2.14)}$$
$i=1,...,r-n$, form a basis of $\Ker \alpha$. We have
$$H_L(\lambda_i,\lambda_j)=\sum_{k=1}^n\bar  z_{ki} z_{kj} - \delta_i^j= \{\bar z^t
z-E_{r-n}\}_{ij}\eqno{(2.15)}$$
hence 2.12 implies 2.5.
\medskip
So, we have constructed a well-defined map from the set of objects of type {\bf (B)} to
the
set of objects of type {\bf (E)}.
\medskip
To construct the inverse map we need a definition. Let $W$ be a $\n C$-vector space. We
denote by $\goth i(W)$ the complex conjugate space together with a map $\goth i: W \to
\goth i(W)$ which is an isomorphism of $\n R$-vector spaces and satisfies $\goth
i(zw)=\bar z \goth i(w)$, $z \in \n C$, $w\in W$.

Let a triple $(L, H_L, \alpha)$ be given. Let $(\Ker \alpha)^\bot\subset L\otimes_K\n C$
be the $H_{L,\n C}$-orthogonal space of $\Ker \alpha$, $\pi_\alpha: L\otimes_K\n C\to
\Ker \alpha$ the projection along $(\Ker \alpha)^\bot$, and let us consider the
composition $\goth i\circ \pi_\alpha: L\otimes_K\n C\to \goth i(\Ker \alpha)$.

We let $V=\n C^n \oplus \goth i(\Ker \alpha)$ (here $\n C^n$ is the target of $\alpha$),
and we define the inclusion $\gamma:L \hookrightarrow V$ by the formula
$$\gamma(l) = (\alpha(l), \goth i\circ \pi_\alpha(l) )$$
The polarization $H$ on $V$ is defined by the same formula (2.10); here $H_L$ is given, so (2.10) defines $\Omega$ for any pair $(\beta(l_1), \beta(l_2))$. Further, $H$ is defined uniquely by $\Omega$ (but clearly we must prove that this $H$ really is a positive hermitian form). The action of $K$ on $V$ is defined by the property that $k*v=kv$ for $v\in \n C^n$ and $k*v=\bar kv$ for $v\in \goth i(\Ker \alpha)$, i.e. $\n C^n=V^+$, $ \goth i(\Ker \alpha)=V^-$.

 Let us prove that the above $V$, $H$, $\gamma(L)$ and the action of $K$ really define an abelian variety up to polarization of type {\bf (B)}. Obviously, the signature of the action of $K$ on $V$ is $(n,r-n)$. Further, for $k\in K$, $l\in L$ we have $\gamma(kl) = k*\gamma(l)$, i.e. $\gamma(L)$ is a $K$-*-lattice. Let $\goth x_1, ...,\goth
x_r$ be a basis of $L$ over $K$ satisfying (2.4). It defines a $n\times (r-n)$-matrix
$z=\{z_{ij}\}$ as follows (this is the same as 2.14):
$$\alpha(\goth x_{n+i})=\sum_{k=1}^n \bar z_{ki}\alpha(\goth x_k)\eqno{(2.16)}$$
(It is easy to prove that $\alpha(\goth x_k)$, $k=1,...,n$ form a basis of $\n C^n$).
Condition (2.5) implies that $z\in \Cal H^3_{n,r-n}$ (calculations coincide with the
ones of 2.15). This means that $\lambda_i$ defined by (2.14) with $z_{**}$ from (2.16) form a basis of Ker $\alpha$. We define
$$\mu_i:=\goth x_{i}\otimes 1 -\sum_{k=1}^{r-n} \goth x_{n+k}\otimes z_{ik}\eqno{(2.17)}$$
$i=1,...,n$. Obviously $H_{L,\n C}(\lambda_i,\mu_j)=0$, this means that $\mu_1,...,\mu_n$ is a basis of $(\Ker \alpha)^\bot$.

Let $Y$ be defined by (2.13) with $z_{**}$ from (2.16). We have $$Y\left(\matrix E_n & -z \\
-\bar z^t&E_{r-n}\endmatrix \right)=\left(\matrix E_n - z\bar z^t& 0 \\ 0&
-\bar z^tz+E_{r-n}\endmatrix \right)$$ hence (2.12) implies that $Y$ is invertible. We define a vector column $e_*$ by the formula $e_*:=Y^{-1}\cdot\gamma(\goth x_*)$. We get immediately using (2.16), (2.17) that $e_1, ... , e_n$ from a basis of $V^+$, $e_{n+1}, ... , e_r$ from a basis of $V^-$. So, results of [Sh] imply that $V$, $H$, $\gamma(L)$ and the action of $K$ really define an abelian variety up to polarization of type {\bf (B)}.
\medskip
The proof of the affirmation that the above two constructions are inverse is practically a tautology. $\square$

\medskip
{\bf Section 3. Exterior powers of abelian varieties with MIQF having $n=1$.}

\nopagebreak
\medskip
Let $(A$, $\iota$, $H)$ be a triple of Theorem 2.6. The associated triple $(L, H_L,
\alpha)$
defines an exact sequence of $\n C$-vector spaces $$0\to \Ker \alpha
\overset{i}\to{\hookrightarrow}
L\otimes_K\n C
\overset{\alpha}\to{\to} \n C^n \to 0$$
Let us take the $k$-th exterior power of this exact sequence (Definition 1.4):
$$0\to \lambda^k(\Ker \alpha) \overset{\lambda^k(i)}\to{\hookrightarrow}
\lambda^k(L)\otimes_K\n C
\overset{\alpha_k}\to{\to} W_k \to 0$$
There exists an Hermitian form $\lambda^k(H_L)$ on $\lambda^k(L)$; recall that
$$\lambda^k(H_L)(l_1\wedge ...\wedge l_k, l'_1\wedge ...\wedge l'_k)=\det \{H_L(l_i\wedge
l'_j)\}$$
It is obvious that if $n=1$ then the restriction of $(-1)^{k}\lambda^k(H_L)$ to
$\lambda^k(\Ker \alpha)$ is positive definite. This means that the triple
$\lambda^k(L)$, $(-1)^{k-1}\lambda^k(H_L)$,
$\alpha_k$ satisfies conditions of Theorem 2.6 and hence defines an abelian variety (up
to isogeny) $\lambda^k(A)$ which is called the $k$-th exterior power of $A$. Its
signature is $r-1 \choose k-1$, $r-1 \choose k$.
\medskip
{\bf Remark 3.1.} It is easy to see that if $n\ne1, r-1$ then this construction cannot be
applied to
abelian varieties of signature $(n,r-n)$; we cannot also define symmetric powers of
abelian varieties with MIQF, as well as their tensor products. This is clearly an analog
of 1.3.
\medskip
{\bf Remark 3.2.} For the function field case there exists a natural definition of the
exterior power of T-motives and 1.7 is a theorem, while for the number
field case there is no such definition hence we must consider the above construction ---
an analog of the Theorem 1.7 --- as the definition of the exterior powers of abelian
varieties with MIQF.
\medskip
{\bf Remark 3.3.} This construction of exterior power can be continued to the corresponding Shimura varieties. Let us recall some definitions of [D]. Attached to a Shimura variety $X$ is an algebraic group $G=G_X$ over $\n Q$ and a map $h=h_X: \n S \to G$ over $\n R$, where $\n S=\Res_{\n C/\n R}(G_m)$, satisfying some conditions. A map $\vf: X \to Y$ of Shimura varieties comes from a map $G_\vf: G_X \to G_Y$ satisfying $G_\vf \circ h_X=h_Y$. If $X$ is a Shimura variety parametrizing abelian varieties with MIQF of signature $(r,s)$ then $G_X=GU(r,s)$. The map $G_\vf: GU(1,r-1)\to GU({r-1 \choose k-1}, {r-1 \choose k})$ corresponding to the above construction is evident. Really, let $L$, $H_L$ be of 2.3 ($n=1$); then $\gamma\in GU(1,r-1)$ is an endomorphism of $L$ (denoted by $\gamma$ as well) satisfying $H_L(\gamma(l_1), \gamma(l_2))=\lambda(\gamma) H_L(l_1,l_2)$. Clearly $G_\vf(\gamma)$ is $\lambda^k(\gamma): \lambda^k(L)\to \lambda^k(L)$. It is an elementary exercise to check all compatibilities.
\medskip

{\bf References}
\nopagebreak
\medskip
[D] Deligne P. Travaux de Shimura. Lect. Notes in Math., 1971,
v.244, p. 123 - 165. Seminaire Bourbaki 1970/71, Expos\'e 389.
\medskip
[G] MR1423131 (97i:11062) Goss, David Basic structures of function
field arithmetic. Springer-Verlag, Berlin, 1996. xiv+422 pp.
\medskip
[GL07] Grishkov, A.; Logachev, D. Duality of Anderson t-motives. https://arxiv.org/pdf/0711.1928.pdf
\medskip
[L] Logachev D., Congruence relations on Anderson T-motives. In preparation.
\medskip
[P] Richard Pink, Hodge structures over function fields. Universit\"at Mannheim.
Preprint. September 17, 1997
\medskip
[Sh] Shimura, Goro On analytic families of polarized abelian
varieties and automorphic functions. Annals of Math., 1 (1963), vol.
78, p. 149 -- 192
\medskip
[W] T. Wedhorn, Congruence relations on some Shimura varieties. J.
reine angew. Math., 524 (2000), p. 43 -- 71
\medskip
\enddocument